\magnification=\magstep1
\input amstex
\documentstyle{amsppt}
\input xy
\xyoption{all}
\voffset=-1.1cm
\define\Z{\Bbb Z}

\topmatter
\title
Rational algebraic K-theory of topological K-theory
\endtitle
\author
Christian Ausoni and John Rognes
\endauthor
\date
August 14th 2007
\enddate
\address
Institute of Mathematics, University of Bonn, DE-53115 Bonn, Germany
\endaddress
\email
ausoni\@math.uni-bonn.de
\endemail
\address
Department of Mathematics, University of Oslo, NO-0316 Oslo, Norway
\endaddress
\email
rognes\@math.uio.no
\endemail
\abstract
We show that after rationalization there is a homotopy fiber sequence
$BBU_\otimes \to K(ku) \to K(\Z)$.  We interpret this as a correspondence
between the virtual $2$-vector bundles over a space $X$ and their
associated anomaly bundles over the free loop space $\Cal L X$.  We also
rationally compute $K(KU)$ by using the
localization sequence, and $K(MU)$ by a method that
applies to all connective $S$-algebras.
\endabstract
\endtopmatter

\define\C{\Bbb C}
\define\Q{\Bbb Q}

\define\im{\operatorname{im}}

\define\Map{\operatorname{Map}}
\define\id{\operatorname{id}}

\define\Spec{\operatorname{Spec}}
\define\Spt{\operatorname{Spt}}

\document

\head
Introduction
\endhead

We are interested in the algebraic $K$-theory $K(ku)$ of the connective
complex $K$-theory spectrum $ku$.  By the calculations of \cite{AR02,
Thm.~0.4} and \cite{Au}, the ``mod~$p$ and~$v_1$'' homotopy of $K(ku)$
is purely $v_2$-periodic, a distinctive homotopy theoretic property it
shares with the spectra representing elliptic cohomology \cite{LRS95}
and topological modular forms \cite{Ho02, \S4}.  The theory of $2$-vector
bundles from \cite{BDR04} and \cite{BDRR} therefore exhibits
$K(ku)$ as a geometrically defined form of elliptic cohomology.
In Section~5 we outline how a $2$-vector bundle with connecting data,
similar to a connection in a vector bundle, is thought to specify
a $1+1$-dimensional conformal field theory.  Since these $2$-vector
bundles are also effective cycles for the form of elliptic cohomology
theory represented by $K(ku)$, we have some justification for referring
to them as elliptic objects, as proposed by Segal \cite{Se89}.

As illustrated by the authors' calculations referred to above, the
arithmetic and homotopy-theoretic information captured by algebraic
$K$-theory becomes more accessible after the introduction of suitable
finite coefficients.  However, for the extraction of $\C$-valued numerical
invariants from a conformal field theory, only the rational homotopy
type of $K(ku)$ will matter.  We are grateful to Ib Madsen and Dennis
Sullivan for insisting that for such geometric applications, we should
first want to compute $K(ku)$ rationally.
To this end we can offer the following theorem.  There is a unit
inclusion map
$$
w \: BBU_\otimes \to BGL_1(ku) \to BGL_\infty(ku)^+ \to K(ku) \,.
$$
Let $\pi \: K(ku) \to K(\Z)$ be induced by the zero-th Postnikov section
$ku \to H\Z$.  The composite $\pi \circ w$ is the constant map to the
base-point of the $1$-component of $K(\Z)$.

\proclaim{Theorem 0.1}
(a)
After rationalization,
$$
BBU_\otimes @>w>> K(ku) @>\pi>> K(\Z)
$$
is a split homotopy fiber sequence.

(b)
The Poincar{\'e} series of $K(ku)$ is
$$
1 + {t^3 \over 1-t^2} + {t^5 \over 1-t^4}
	= 1 + {t^3 + 2t^5 \over 1-t^4} \,.
$$

(c)
There is a rational determinant map
$$
{\det}_\Q \: BGL_\infty(ku)^+ \to BGL_1(ku)_\Q
$$
that, in its relative form for $ku \to H\Z$, rationally splits $w$.
\endproclaim

By the Poincar{\'e} series of a space $X$ of finite type, we mean the formal
power series $\sum_{n\ge0} r_n t^n$ in $\Z[[t]]$, where $r_n$ is the
rank of $\pi_n(X)$.
Theorem~0.1 is proved by assembling Theorem~2.5(a) and Theorem~4.8(a).
The other parts of those theorems prove similar results for $K(ko)$ and
$K(\ell)$, where $\ell = BP\langle1\rangle$.

The splitting of $K(ku)_\Q$ shows that for a (virtual) $2$-vector bundle
over $X$, represented by a map $\Cal E \: X \to K(ku)$, the rational
information splits into two pieces.  The less interesting piece is the
decategorified information carried by the dimension bundle $\dim(\Cal E)
= \pi \circ \Cal E \: X \to K(\Z)$.  The more interesting piece is the
determinant bundle
$$
|\Cal E| = {\det}_\Q \circ \Cal E \: X \to (BBU_\otimes)_\Q \,.
$$
To specify $|\Cal E|$ is equivalent to specifying a rational virtual
vector bundle $\Cal H \: \Cal L X \to (BU_\otimes)_\Q$ over the free
loop space $\Cal L X = \Map(S^1, X)$, called the ``anomaly bundle'',
subject to a coherence condition relating the composition of free
loops, when defined, to the tensor product of virtual vector spaces.
See diagram~(5.3).

The conclusion is that for rational purposes the information in
a $2$-vector bundle~$\Cal E$ over $X$ is the same as that in its
anomaly bundle $\Cal H$ over $\Cal L X$ (subject to the indicated
coherence condition, which we think of as implicit, and together with
the dimension bundle $\dim(\Cal E)$ over $X$, which we tend to ignore).
In physical language, the fiber of the anomaly bundle at a free loop
$\gamma \: S^1 \to X$ plays the role of the state space of $\gamma$
viewed as a closed string in $X$.  The advantage of $2$-vector bundles
over their homotopy-theoretic alternatives, such as representing maps
to classifying spaces or bundles of $ku$-modules, is that they are
geometrically modeled in terms of vector bundles, rather than virtual
vector bundles.  This seems to become an essential virtue when one
wants to treat differential-geometric structures like connections on
these bundles.

We also compute the rational algebraic $K$-theory $K(KU)$ of the
periodic complex $K$-theory spectrum $KU$. To this end we evaluate the
(rationalized) transfer map $\pi_*$ in the localization sequence
$$
K(\Z) @>\pi_*>> K(ku) @>\rho>> K(KU)
$$
predicted by the second author and established by Blumberg
and Mandell \cite{BM}.

\proclaim{Theorem 0.2}
(a)
There is a rationally split homotopy fiber sequence of infinite loop
spaces
$$
K(ku) @>\rho>> K(KU) @>\partial>> BK(\Z)
$$
where $\rho$ is induced by the connective cover map $ku \to KU$,
and $BK(\Z)$ denotes the first connected delooping of $K(\Z)$.

(b)
The Poincar{\'e} series of $K(KU)$ is
$$
(1+t) + {t^3 + 2t^5 + t^6 \over 1-t^4} \,.
$$
\endproclaim

See Theorem~2.12 for our proof.

As stated, Theorems~0.1 and~0.2 only concern the algebraic $K$-theory
of topological $K$-theory, but we develop our proofs in the greater
generality of arbitrary connective $S$-algebras.  In Section~1 we
observe how the calculation by Goodwillie \cite{Go86} of the relative
rational algebraic $K$-theory for a $1$-connected map $R \to \pi_0 R$
of simplicial rings (which generalizes earlier calculations by Hsiang
and Staffeldt \cite{HS82} for simplicial group rings), also applies to
determine the relative rational algebraic $K$-theory for a $1$-connected
map $A \to H\pi_0 A$ of connective $S$-algebras.  The answer is given
in terms of negative cyclic homology; see Theorem~1.5 and Corollary~1.6.

When $\pi_0 A$ is close to $\Z$, and $A \to H\pi_0 A$ is a ``rational
de\,Rham equivalence'', we get a very simple expression for the relative
rational algebraic $K$-theory as the image of Connes' $B$-operator on
Hochschild homology; see Proposition~1.8 and Corollary~1.9.  These
hypotheses apply to a number of interesting examples of connective
$S$-algebras, including the $K$-theory spectra $ku$, $ko$ and $\ell$,
and the bordism spectra $MU$, $MSO$ and $MSp$.  We work these examples
out in Theorem~2.5 and Theorem~3.4, respectively.  

In Section~4 we consider the unit inclusion map $w \: BGL_1(A) \to K(A)$.
For commutative $A$, the rationalization $A_\Q$ is equivalent as a
commutative $H\Q$-algebra to the Eilenberg--Mac\,Lane spectrum $HR$
of a commutative simplicial $\Q$-algebra, so we can use the determinant
$GL_n(R) \to GL_1(R)$ to define a rational determinant map
$$
{\det}_\Q \: BGL_\infty(A)^+ \to BGL_1(A)_\Q \,.
$$
We show in Proposition~4.7 that the composite $\det_\Q \circ w$ is the
rationalization map, and apply this in Theorem~4.8 to show that $w$
induces a rational equivalence from $BBU_\otimes$ to the homotopy fiber
of $\pi \: K(ku) \to K(\Z)$, and similarly for $ko$ and $\ell$.
This last step is a counting argument; it does not apply for $MU$
or the other bordism spectra.

\subhead Acknowledgments \endsubhead
In an earlier version of this paper, we emphasized a trace map to
$THH(ku)$ over the determinant map to $(BBU_\otimes)_\Q$, in order
to detect the image of $w$ in $K(ku)$.  We are grateful to Bj{\o}rn
Dundas for reminding us of the existence of determinants for commutative
simplicial rings, which is half of the basis for the existence of the
map $\det_\Q$ defined in Lemma~4.6.
We are also grateful to Mike Mandell and Brooke Shipley for help with
some of the references concerning commutative simplicial $\Q$-algebras
given in Subsection~1.1.

\head
\S1. Rational algebraic $K$-theory of connective $S$-algebras
\endhead

\subhead 1.1. $S$-algebras \endsubhead
We work in one of the modern symmetric monoidal categories of spectra
\cite{EKMM97}, \cite{Ly99}, \cite{HSS00}, \cite{MMSS01}, which we shall
refer to as $S$-modules.  The monoids (resp.~commutative monoids) in this
category are called $S$-algebras (resp.~commutative $S$-algebras), and are
equivalent to the $A_\infty$ ring spectra (resp.~$E_\infty$ ring spectra)
considered since the 1970's.  The Eilenberg--Mac\,Lane functor $H \: R \mapsto
HR$ maps the category of simplicial rings (resp.~commutative simplicial
rings) to the category of $S$-algebras (resp.~commutative $S$-algebras).

Schwede proved in \cite{Schw99, 4.5} that $H$ is part of a Quillen
equivalence from the category of simplicial rings to the category of
connective $H\Z$-algebras.  There is a similar equivalence between the
category of commutative simplicial $\Q$-algebras and the category of
connective commutative $H\Q$-algebras.

One form of the latter equivalence appears in \cite{KM95, II.1.3}.  In a
little more detail, the category of connective commutative $H\Q$-algebras
is ``connective Quillen equivalent'' \cite{MMSS01, p.~445} to the category
of connective $E_\infty$ $H\Q$-ring spectra \cite{EKMM97, II.4}, and
connective $E_\infty$ $H\Q$-ring spectra are the $E_\infty$ objects in
connective $H\Q$-modules, which are Quillen equivalent to $E_\infty$
simplicial $\Q$-algebras \cite{Schw99, 4.4}.  The monads defining
$E_\infty$ algebras and commutative algebras in simplicial $\Q$-modules
are weakly equivalent, since for every $j\ge0$ the group homology of
$\Sigma_j$ with coefficients in any $\Q$-module is concentrated in
degree zero.  Hence $E_\infty$ simplicial $\Q$-algebras are Quillen
equivalent to commutative simplicial $\Q$-algebras \cite{Ma03, 6.7}.

The homotopy categories of commutative simplicial rings and connective
commutative $H\Z$-algebras are not equivalent.

\subhead 1.2. Linearization \endsubhead
Let $A$ be a connective $S$-algebra.  We write $\pi = \pi_A \: A \to
H\pi_0 A$ for its zero-th Postnikov section, and define the linearization
map $\lambda = \lambda_A \: A \to H\Z \wedge A$ to be $\pi \wedge id \:
A \cong S \wedge A \to H\Z \wedge A$.  It is a $\pi_0$-isomorphism and a
rational equivalence of connective $S$-algebras.  For each (simplicial or
topological) monoid $G$ let $S[G] = \Sigma^\infty G_+$ be its unreduced
suspension spectrum.  For $A = S[G]$, the linearization map $\lambda \:
S[G] \to H\Z \wedge S[G] \cong H\Z[G]$ agrees with the map considered
by Waldhausen \cite{Wa78, p.~43}.

In general, $H\Z \wedge A$ is a connective $H\Z$-algebra, so by the first
Quillen equivalence above there is a naturally associated simplicial
ring $R$ with $H\Z \wedge A \simeq HR$.  For connective commutative $A$,
the rationalization $A_\Q = H\Q \wedge A$ is a connective commutative
$H\Q$-algebra, so by the second Quillen equivalence above there is a
naturally associated commutative simplicial $\Q$-algebra $R$ with $A_\Q
\simeq HR$.

\subhead 1.3. Algebraic $K$-theory \endsubhead
For a general $S$-algebra $A$, the algebraic $K$-theory space $K(A)$ can
be defined as $\Omega |h S_\bullet \Cal C_A|$, where $\Cal C_A$ is the
category of finite cell $A$-module spectra and their retracts, $S_\bullet$
denotes Waldhausen's $S_\bullet$-construction \cite{Wa85, \S1.3}, and
$|h(-)|$ indicates the nerve of the subcategory of weak equivalences.
By iterating the $S_\bullet$-construction, we may also view $K(A)$
as a spectrum.  For connective $S$-algebras, $K(A)$ can alternatively
be defined in terms of Quillen's plus-construction as $K_0(\pi_0 A)
\times BGL_\infty(A)^+$, and then the two definitions are equivalent,
see \cite{EKMM97, VI.7.1}.  We write $K(R)$ for $K(HR)$, and similarly
for other functors.

Any map $A \to A'$ of connective $S$-algebras that is a
$\pi_0$-isomorphism and a rational equivalence induces a rational
equivalence $K(A) \to K(A')$, by \cite{Wa78, 2.2}.  The proof goes by
observing that $BGL_n(A) \to BGL_n(A')$ is a rational equivalence
for each $n$.  In particular, for $R$ with $H\Z \wedge A \simeq HR$
there is a natural rational equivalence $\lambda \: K(A) \to K(H\Z
\wedge A) \simeq K(R)$.  For $X \simeq BG$, Waldhausen writes $A(X)$ for
$K(S[G])$, and $\lambda \: A(X) \to K(\Z[G])$ is a rational equivalence.
In this case, $A(X)$ can also be defined as the algebraic $K$-theory
of a category $\Cal R_f(X)$ of suitably finite retractive spaces over $X$,
see \cite{Wa85, \S2.1}.

\subhead 1.4. Cyclic homology \endsubhead
There is a natural trace map $tr \: K(A) \to THH(A)$ to the topological
Hochschild homology of $A$, see \cite{BHM93, \S3}.  The target is a cyclic
object in the sense of Connes, hence carries a natural $S^1$-action.
There exists a model for the trace map that factors through the fixed
points of this circle action \cite{Du04}, hence it also factors through
the homotopy fixed points $THH(A)^{hS^1}$.  We get a natural commutative
triangle
$$
\xymatrix{
K(A) \ar[r]^-\alpha \ar[dr]_{tr} & THH(A)^{hS^1} \ar[d]^F \\
& THH(A)
}
$$
where the Frobenius map $F$ forgets about $S^1$-homotopy invariance.  For
any simplicial ring $R$ there are natural isomorphisms $THH_*(HR_\Q) \cong
HH_*(R \otimes \Q)$ (Hoch\-schild homology) and $THH_*(HR_\Q)^{hS^1} \cong
HC^-_*(R \otimes \Q)$ (negative cyclic homology).  See e.g.~\cite{EKMM97,
IX.1.7} and \cite{CJ90, 1.3(3)}.  With these identifications, the
triangle above realizes the commutative diagram of \cite{Go86, II.3.1}.
In \cite{Go86, II.3.4}, Goodwillie proved:

\proclaim{Theorem 1.5}
Let $f \: R \to R'$ be a map of simplicial rings, with $\pi_0 R \to
\pi_0 R'$ a surjection with nilpotent kernel.  Then
$$
\xymatrix{
K(R)_\Q \ar[r]^-{\alpha} \ar[d]_f & HC^-(R \otimes \Q) \ar[d]^f \\
K(R')_\Q \ar[r]^-{\alpha} & HC^-(R' \otimes \Q)
}
$$
is homotopy Cartesian, i.e., the map of vertical homotopy fibers
$$
\alpha \: K(f)_\Q \to HC^-(f \otimes \Q)
$$
is an equivalence.
\endproclaim

Here we write $K(f)$ for the homotopy fiber of $K(R) \to K(R')$,
so that there is a long exact sequence
$$
\dots \to K_{*+1}(R') \to K_*(f) \to K_*(R) \to K_*(R') \to \dots \,,
$$
and similarly for other functors from ($S$-)algebras to spaces.
(Goodwillie writes $K(f)$ for a delooping of our $K(f)$, but we
need to emphasize fibers over cofibers.)  We write $K(R)_\Q$ for
the rationalization of $K(R)$, and similarly for other spaces and
$S$-algebras.

\proclaim{Corollary 1.6}
Let $g \: A \to A'$ be a map of connective $S$-algebras, with
$\pi_0 A \to \pi_0 A'$ a surjection with nilpotent kernel,  Then
$$
\xymatrix{
K(A)_\Q \ar[r]^-{\alpha} \ar[d]_g & THH(A_\Q)^{hS^1} \ar[d]^g \\
K(A')_\Q \ar[r]^-{\alpha} & THH(A'_\Q)^{hS^1}
}
$$
is homotopy Cartesian, i.e., the map of vertical homotopy fibers
$$
\alpha \: K(g)_\Q \to THH(g_\Q)^{hS^1}
$$
is an equivalence.
\endproclaim

\demo{Proof}
Given $g \: A \to A'$ we find $f \: R \to R'$ with $H\Z \wedge A \simeq HR$
and $H\Z \wedge A' \simeq HR'$ making the diagram
$$
\xymatrix{
A \ar[r]^\lambda \ar[d]_g & HR \ar[d]^{Hf} \\
A' \ar[r]^\lambda & HR'
}
$$
homotopy commute.  Then $\lambda \: K(A) \to K(R)$ is a rational
equivalence and $A_\Q \simeq H(R \otimes \Q)$, so the square in the
corollary is equivalent to the square in Goodwillie's theorem.
\qed
\enddemo

\subhead 1.7. De\,Rham homology \endsubhead
The (spectrum level) circle action on $THH(A)$ induces a suspension
operator $d \: THH_*(A) \to THH_{*+1}(A)$, analogous to Connes' operator
$B \: HH_*(R) \to HH_{*+1}(R)$.  When $A_\Q \simeq H(R \otimes \Q)$,
these operators are compatible under the isomorphism $THH_*(A_\Q) \cong
HH_*(R \otimes \Q)$.  In general $dd = d\eta$ is not zero \cite{He97,
1.4.4}, where $\eta$ is the stable Hopf map, but in the algebraic case
$BB = 0$, so one can define the de\,Rham homology
$$
H^{dR}_*(R) = \ker(B)/\im(B)
$$
of a simplicial ring $R$ as the homology of $HH_*(R)$ with respect to
the $B$-operator.

For a map $g \: A \to A'$ of $S$-algebras, the homotopy fiber $THH(g)$ of
$THH(A) \to THH(A')$ inherits a circle action and associated suspension
operator.  Similarly, for a map $f \: R \to R'$ of simplicial rings
there is a relative $B$-operator acting on the term $HH_*(f)$ in
the long exact sequence
$$
\dots \to HH_{*+1}(R') \to HH_*(f) \to HH_*(R) \to HH_*(R') \to \dots \,,
$$
and we define $H^{dR}_*(f)$ to be the homology of $HH_*(f)$ with
respect to this $B$-operator.  We say that $f \: R \to R'$ is a de\,Rham
equivalence if $H^{dR}_*(f) = 0$, and that $f$ is a rational de\,Rham
equivalence if $H^{dR}_*(f \otimes \Q) = 0$.  If we assume that $HH_*(R)
\to HH_*(R')$ is surjective in each degree, then there is a long exact
sequence
$$
\dots \to H^{dR}_{*+1}(R') \to H^{dR}_*(f) \to H^{dR}_*(R) \to
H^{dR}_*(R') \to \dots \,,
$$
in which case $f$ is a de\,Rham equivalence if and only if $H^{dR}_*(R)
\to H^{dR}_*(R')$ is an isomorphism in every degree.

\proclaim{Proposition 1.8}
If $f \: R \to R'$ is a de\,Rham equivalence,
then there is an exact sequence
$$
0 \to HC^-_*(f) @>F>> HH_*(f) @>B>> HH_{*+1}(f)
$$
that identifies $HC^-_*(f)$ with $\ker(B) \subset HH_*(f)$.
\endproclaim

\demo{Proof}
By analogy with the homotopy fixed point spectral sequence for
$THH(g)^{hS^1}$, there is a second quadrant homological spectral sequence
$$
E^2_{**} = \Q[t] \otimes HH_*(f) \Longrightarrow HC^-_*(f)
$$
with $t \in E^2_{-2,0}$ and $d^2(t^i \cdot x) = t^{i+1} \cdot B(x)$
for all $x \in HH_*(f)$, $i\ge0$.  So $E^3_{**}$ is the sum of $\ker(B)
\subset HH_*(f)$ in the zero-th column and a copy of $H^{dR}_*(f)$ in each
even, negative column.  By assumption the latter groups are all zero, so
the spectral sequence collapses to the zero-th column at the $E^3$-term.
The Frobenius $F$ is the edge homomorphism for this spectral sequence,
and the assertion follows.
\qed
\enddemo

\proclaim{Corollary 1.9}
Let $A$ be a connective $S$-algebra such that $\pi_0 A$ is any
localization of the integers, and let $R$ be a simplicial $\Q$-algebra
with $A_\Q \simeq HR$.

(a)~The homotopy fiber sequence
$$
K(\pi_A) \to K(A) @>\pi_A>> K(\pi_0 A)
$$
is rationally split, where $\pi_A \: A \to H\pi_0 A$ is the zero-th
Postnikov section.

(b)~There are equivalences
$$
K(\pi_A)_\Q @>\alpha>\simeq> THH(\pi_{A\Q})^{hS^1} \simeq HC^-(\pi_R)
\,,
$$
where $\pi_R \: R \to \pi_0 R = \Q$ is the zero-th Postnikov section.

Suppose furthermore that $H^{dR}_*(R) \cong \Q$ is trivial in positive
degrees.

(c)~The map $\pi_R$ is a de\,Rham equivalence, and the Frobenius map
identifies $HC^-_*(\pi_R)$ with the positive-degree part of
$$
\ker(B) \subset HH_*(R) \cong THH_*(A) \otimes \Q \,.
$$
That part is also equal to $\im(B) \subset THH_*(A) \otimes \Q$.

(d)~The trace map $tr \: K(A) \to THH(A)$ induces the composite
identification of $K_*(\pi_A) \otimes \Q$ with the positive-degree part
of $\ker(B) \subset THH_*(A) \otimes \Q$.
\endproclaim

\demo{Proof}
(a) Write $\pi_0 A = \Z_{(P)}$ for some (possibly empty) set of primes
$P$.  The unit map $i \: S \to A$ factors through $S_{(P)}$, and the
composite map $S_{(P)} \to A \to H\pi_0 A$ is a $\pi_0$-isomorphism and
a rational equivalence.  Hence the composite
$$
K(S_{(P)}) \to K(A) @>\pi_A>> K(\pi_0 A)
$$
is a rational equivalence.

(b) The map $\pi_A \: A \to H\pi_0 A$ induces the identity on $\pi_0$, 
so $\alpha$ is an equivalence by Corollary~1.6.  We recalled the
second identification in Subsection~1.4.  It is clear that
$\pi_0 R = \pi_0 A_\Q = \pi_0 A \otimes \Q = \Q$.

(c) Since $HH_*(\Q) = \Q$ is trivial in positive degrees, the map
$HH_*(R) \to HH_*(\Q)$ is surjective in each degree, so $\pi_R$ is a
de\,Rham equivalence if (and only if) $H^{dR}_*(R) \cong H^{dR}_*(\Q)
= \Q$ is trivial in all positive degrees.  The homotopy fiber sequence
$$
HH(\pi_R) \to HH(R) \to HH(\Q)
$$
identifies $HH_*(\pi_R)$ with the positive-degree part of $HH_*(R)$,
so $\ker(B) \subset HH_*(\pi_R)$ is the  positive-degree part of
$\ker(B) \subset HH_*(R)$.
The identification $\im(B) = \ker(B)$ in positive degrees is of
course equivalent to the vanishing of $H^{dR}_*(R)$ in positive
degrees.

(d) The trace map factors as $tr = F \circ \alpha$.
\qed
\enddemo

\head \S2. Examples from topological $K$-theory \endhead

\subhead 2.1. Connective $K$-theory spectra \endsubhead
Let $ku$ be the connective complex $K$-theory spectrum, $ko$ the
connective real $K$-theory spectrum, and $\ell = BP\langle 1\rangle$
the Adams summand of $ku_{(p)}$, for $p$ an odd prime.  These are all
commutative $S$-algebras.  We write $\Omega^\infty ku = BU \times \Z$,
$\Omega^\infty ko = BO \times \Z$ and $\Omega^\infty \ell = W \times
\Z_{(p)}$ for the underlying infinite loop spaces (see \cite{Ma77, V.3-4}).
The homotopy units form infinite loop spaces, namely $GL_1(ku) =
BU_\otimes \times \{\pm1\}$, $GL_1(ko) = BO_\otimes \times \{\pm1\}$
and $GL_1(\ell) = W_\otimes \times \Z_{(p)}^\times$.  The homotopy
algebras are $\pi_* ku = \Z[u]$ with $|u|=2$, $\pi_* ko = \Z[\eta,
\alpha, \beta]/(2\eta, \eta^3, \eta\alpha, \alpha^2 {-} 4\beta)$ with
$|\eta|=1$, $|\alpha|=4$, $|\beta|=8$, and $\pi_* \ell = \Z_{(p)}[v_1]$
with $|v_1| = 2p-2$.  The complexification map $ko \to ku$ takes $\eta$
to $0$, $\alpha$ to $2u^2$ and $\beta$ to $u^4$.  The inclusion $\ell
\to ku_{(p)}$ takes $v_1$ to $u^{p-1}$.

\proclaim{Proposition 2.2}
(a) There are $\pi_0$-isomorphisms and rational equivalences
$$
\kappa \: S[\Omega S^3] \to S[K(\Z,2)] \to ku
$$
of $S$-algebras, so $ku_\Q \simeq H\Q[\Omega S^3]$ as homotopy commutative
$H\Q$-algebras, where $\Q[\Omega S^3]$ is a simplicial $\Q$-algebra,
and $ku_\Q \simeq H\Q[K(\Z,2)]$ as commutative $H\Q$-algebras, where
$\Q[K(\Z,2)]$ is a commutative simplicial $\Q$-algebra.

(b) There are $\pi_0$-isomorphisms and rational equivalences $\bar \alpha
\: S[\Omega S^5] \to ko$ and $\bar v_1 \: S_{(p)}[\Omega S^{2p-1}] \to
\ell$, so $ko_\Q \simeq H\Q[\Omega S^5]$ and $\ell_\Q \simeq H\Q[\Omega
S^{2p-1}]$ as $H\Q$-algebras, where $\Q[\Omega S^5]$ and $\Q[\Omega
S^{2p-1}]$ are simplicial $\Q$-algebras.

(c) In particular, there is a rational equivalence $\kappa \: A(S^3)
\to K(ku)$ of $S$-algebras, a rational equivalence $A(K(\Z,3)) \to K(ku)$
of commutative $S$-algebras, and a rational equivalence $\bar \alpha \:
A(S^5) \to K(ko)$ of $S$-modules.
\endproclaim

\demo{Proof}
(a) Let $BS^3 \to K(\Z,4)$ represent a generator of $H^4(BS^3)$.
It induces a double loop map $\Omega S^3 \to K(\Z,2)$, such that the
composite $S^2 \to \Omega S^3 \to K(\Z,2)$ represents a generator of
$\pi_2 K(\Z,2)$.  The inclusions $K(\Z,2) \simeq BU(1) \to BU_\otimes \to
GL_1(ku)$ are infinite loop maps, and the generator of $\pi_2 K(\Z,2)$
maps to a generator of $\pi_2 GL_1(ku)$.  By adjunction we have an $E_2$
ring spectrum map $S[\Omega S^3] \to S[K(\Z,2)]$ and an $E_\infty$ ring
spectrum map $S[K(\Z,2)] \to ku$, with composite the $E_2$ ring spectrum
map $\kappa \: S[\Omega S^3] \to ku$.

These are rational equivalences, because $\pi_* S[\Omega S^3] \otimes \Q
\cong H_*(\Omega S^3; \Q) \cong \Q[x]$, $H_*(K(\Z,2); \Q) \cong \Q[b]$
and $\pi_* ku \otimes \Q = \Q[u]$, with $\kappa$ mapping $x$ via $b$
to $u$.  We may take the Kan loop group of $S^3$ (a simplicial group, see
e.g.~\cite{Wa96}) as our model for $\Omega S^3$, and rigidify $\kappa$
to a map of $S$-algebras.  Following \cite{FV}, there remains an $E_1
= A_\infty$ operad action on these $S$-algebras and $\kappa$, which
in particular implies that $\kappa \: A(S^3) \to K(ku)$ is homotopy
commutative.

(b) For the real case, let $S^4 \to BO_\otimes \subset GL_1(ko)$ represent
a generator of $\pi_4 GL_1(ko)$.  By the loop structure on the target,
it extends to a loop map $\Omega S^5 \to GL_1(ko)$, with left adjoint an
$A_\infty$ ring spectrum map $\bar \alpha \: S[\Omega S^5] \to ko$.  It is
a rational equivalence, because $\pi_* S[\Omega S^5] \otimes \Q \cong
H_*(\Omega S^5; \Q) \cong \Q[y]$ and $\pi_* ko \otimes \Q = \Q[\alpha]$,
with $\bar\alpha$ mapping $y$ to $\alpha$.  We interpret $\Omega S^5$ as
the Kan loop group, and form the simplicial $\Q$-algebra $\Q[\Omega S^5]$
as its rational group ring.  The Adams summand case is entirely similar,
starting with a map $S^{2p-2} \to W_\otimes \subset GL_1(\ell)$.

(c) By \cite{FV} and naturality there is an $A_\infty$ operad action
on the induced map of spectra $A(S^3) = K(S[\Omega S^3]) \to K(ku)$
(rather than of spaces), which we can rigidify to a map of $S$-algebras.
The $S$-algebra multiplication $A(S^3) \wedge A(S^3) \to A(S^3)$ is
induced by the group multiplication $S^3 \times S^3 \to S^3$.
\qed
\enddemo

\proclaim{Lemma 2.3}
(a) For any integer $n\ge1$ the simplicial $\Q$-algebra $R = \Q[\Omega
S^{2n+1}]$ has Hochschild homology
$$
HH_*(R) \cong \Q[x] \otimes E(dx)
$$
with $|x|=2n$, where Connes' $B$-operator satisfies $B(x) = dx$.
Here $E(-)$ denotes the exterior algebra.
(b) The de\,Rham homology $H^{dR}_*(R) \cong \Q$ is concentrated
in degree zero, so $\pi_R \: R \to \Q$ is a de\,Rham equivalence.
(c) The positive-degree part of $\ker(B) \subset HH_*(R)$ is
$$
\Q[x] \{dx\} = \Q \{ dx, x \, dx, x^2 \, dx, \dots \} \,.
$$
\endproclaim

\demo{Proof}
(a) The Hochschild filtration on the bisimplicial $\Q$-algebra $HH(R)$ yields
a spectral sequence
$$
E^2_{**} = HH_*(\pi_*(R)) \Longrightarrow HH_*(R) \,,
\tag 2.4
$$
and $\pi_*(R) = \Q[x]$ with $|x|=2n$.  The Hochschild homology of
this graded commutative ring is $\Q[x] \otimes E(dx)$, where $dx
\in E^2_{1,2n}$ is the image of $x$ under Connes' $B$-operator.
The spectral sequence collapses at that stage, for bidegree reasons.

(b,c) The $B$-operator is a derivation, hence takes $x^m$ to $m x^{m-1}
\, dx$ for all $m\ge0$.  It follows easily that the de\,Rham homology
is trivial in positive degrees, and that $\ker B$ is as indicated.
\qed
\enddemo

By combining Corollary~1.9, Proposition~2.2 and Lemma~2.3, we obtain
the following result.

\proclaim{Theorem 2.5}
(a)
There is a rationally split homotopy fiber sequence
$$
K(\pi_{ku}) \to K(ku) @>\pi>> K(\Z)
$$
and the trace map $tr \: K(ku) \to THH(ku)$ identifies
$$
K_*(\pi_{ku}) \otimes \Q \cong \Q[u] \{du\}
$$
with its image in $THH_*(ku) \otimes \Q \cong \Q[u] \otimes E(du)$.
Here $|u|=2$ and $|du|=3$, so $K(\pi_{ku})$ has Poincar{\'e} series
$t^3/(1-t^2)$.

(b)
Similarly, there are rationally split homotopy fiber sequences
$$
\align
K(\pi_{ko}) &\to K(ko) @>\pi>> K(\Z) \\
K(\pi_\ell) &\to K(\ell) @>\pi>> K(\Z_{(p)})
\endalign
$$
and the trace maps identify
$$
\align
K_*(\pi_{ko}) \otimes \Q &\cong \Q[\alpha] \{d\alpha\} \\
K_*(\pi_\ell) \otimes \Q &\cong \Q[v_1] \{dv_1\}
\endalign
$$
with their images in $THH_*(ko) \otimes \Q \cong \Q[\alpha] \otimes
E(d\alpha)$ and $THH_*(\ell) \otimes \Q \cong \Q[v_1] \otimes E(dv_1)$,
respectively.  Hence $K(\pi_{ko})$ has Poincar{\'e} series $t^5/(1-t^4)$,
whereas $K(\pi_\ell)$ has Poincar{\'e} series $t^{2p-1}/(1-t^{2p-2})$.
\endproclaim

\remark{Remark 2.6}
The Poincar{\'e} series of $K(\Z)$ is $1 + t^5/(1-t^4)$ by Borel's
calculation \cite{Bo74}.  Hence the (common) Poincar{\'e} series of $K(ku)$
and $A(S^3)$ is
$$
1 + t^3/(1-t^2) + t^5/(1-t^4) = 1 + (t^3+2t^5)/(1-t^4) \,,
$$
whereas the Poincar{\'e} series of $K(ko)$ and $A(S^5)$ is $1 +
2t^5/(1-t^4)$.  More generally, we recover the Poincar{\'e} series
$1 + t^5/(1-t^4) + t^{2n+1}/(1-t^{2n})$ of $A(S^{2n+1})$ for $n\ge1$,
from \cite{HS82, Cor.~1.2}.  The group $K_1(\Z_{(p)})$ is not finitely
generated, so we do not discuss the Poincar{\'e} series of $K(\ell)$.
\endremark

\subhead 2.7. Periodic $K$-theory spectra \endsubhead
Let $KU$ be the periodic complex $K$-theory spectrum, $KO$ the periodic
real $K$-theory spectrum, and $L = E(1)$ the Adams summand of $KU_{(p)}$,
for $p$ an odd prime.  We have maps of commutative $S$-algebras
$$
H\Z @<\pi<< ku @>\rho>> KU
$$
with associated maps of ``brave new'' affine schemes \cite{TV, \S2}
$$
\Spec(\Z) @>\pi>> \Spec(ku) @<\rho<< \Spec(KU) \,.
\tag 2.8
$$
Let $i \: S \to S[\Omega S^3]$ and $c \: S[\Omega S^3] \to S$ be induced
by the inclusion map $* \to S^3$ and the collapse map $S^3 \to *$,
respectively.  We have a map of horizontal cofiber sequences
$$
\xymatrix{
\Sigma^2 S[\Omega S^3] \ar[r]^-x \ar[d]^{\Sigma^2 \kappa} &
S[\Omega S^3] \ar[r]^-c \ar[d]^\kappa &
S \ar[d]^\lambda \\
\Sigma^2 ku \ar[r]^-u & ku \ar[r]^-\pi & H\Z
}
\tag 2.9
$$
where the top row exhibits $S$ as a two-cell $S[\Omega S^3]$-module,
and the bottom row exhibits $H\Z$ as a two-cell $ku$-module.
(In each case, the two cells are in dimension zero and three.)  There are
algebraic $K$-theory transfer maps $c_* \: A(*) \to A(S^3)$ and $\pi_*
\: K(\Z) \to K(ku)$ (with a lower star, in accordance with the variance
conventions from algebraic geometry and~(2.8)), that are induced by
the functors that view finite cell $S$-modules as finite cell $S[\Omega
S^3]$-modules, and finite cell $H\Z$-modules as finite cell $ku$-modules,
respectively.  In terms of retractive spaces, $c_*$ is
induced by the exact functor $\Cal R_f(*) \to \Cal R_f(S^3)$ that takes a
pointed space $X \rightleftarrows *$ to the retractive space
$X \times S^3 \rightleftarrows S^3$.  The transfer maps are compatible,
by~(2.9), so we have a commutative diagram with vertical rational
equivalences
$$
\xymatrix{
A(*) \ar[r]^{c_*} \ar[d]^\lambda &
A(S^3) \ar[d]^\kappa \\
K(\Z) \ar[r]^{\pi_*} & K(ku) \ar[r]^\rho & K(KU) \,.
}
\tag 2.10
$$
The bottom row is a homotopy fiber sequence by the
localization theorem of \cite{BM}.

\proclaim{Lemma 2.11}
The transfer map $c_* \: A(*) \to A(S^3)$ 
is null-homotopic, as a map of $A(*)$-module spectra.
The transfer map $\pi_* \: K(\Z) \to K(ku)$
is rationally null-homotopic, again as a map of $A(*)$-module spectra.
\endproclaim

\demo{Proof}
The projection formula asserts that $c_*$ is an $A(S^3)$-module
map, where $c \: A(S^3) \to A(*)$ makes $A(*)$ an $A(S^3)$-module.
Restricting the module structures along $i \: A(*) \to A(S^3)$, we see
that $c_*$ is a map of $A(*)$-module spectra, and the source is a free
$A(*)$-module of rank one.  Hence it suffices to show that $c_*$ takes a
generator of $\pi_0 A(*)$, represented say by $S^0 \rightleftarrows *$,
to zero in $\pi_0 A(S^3) \cong \Z$.  But $c_*$ maps that generator to
the class of $S^0 \times S^3 \rightleftarrows S^3$, which corresponds
to its Euler characteristic $\chi(S^3) = 0$.

The conclusion for $\pi_*$ follows from that for $c_*$, via the rational
equivalences $\lambda$ and $\kappa$.
\qed
\enddemo

Note the utility of the comparison with $A$-theory at this point, since
we do not have an $S$-algebra map $K(\Z) \to K(ku)$ that is analogous
to $i \: A(*) \to A(S^3)$.

\proclaim{Theorem 2.12}
There are rationally split homotopy fiber sequences
$$
\align
K(ku) &@>\rho>> K(KU) @>\partial>> BK(\Z) \\
K(\ell) &@>\rho>> K(L) @>\partial>> BK(\Z_{(p)})
\endalign
$$
of infinite loop spaces.
Hence the Poincar{\'e} series of $K(KU)$ is
$$
(1+t) + (t^3 +2t^5 + t^6)/(1-t^4) \,.
$$
\endproclaim

\demo{Proof}
The claims for $KU$ follow by combining Theorem~2.5(a) and Lemma~2.11.
The proof of the claim for $L$ is completely similar, using that
$H\Z_{(p)}$ is a two-cell $\ell$-module, with cells in dimension zero and
$(2p-1)$.  By \cite{BM} there is a homotopy fiber sequence $K(\Z_{(p)})
\to K(\ell) \to K(L)$.
\qed
\enddemo

\remark{Remark 2.13}
We do not know how to relate $K(ko)$ with $K(KO)$, so we do not have a
rational calculation of $K(KO)$.  However, $KO \to KU$ is a $\Z/2$-Galois
extension of commutative $S$-algebras, in the sense of \cite{Ro,
\S4.1}, so it is plausible that $K(KO) \to K(KU)^{h\Z/2}$ is close to
an equivalence.  Here $\Z/2$ acts on $KU$ by complex conjugation, and
$$
\pi_*(K(KU)^{h\Z/2}) \otimes \Q \cong [ K_*(KU) \otimes \Q ]^{\Z/2} \,.
$$
The conjugation action on $ku$ fixes $K(\Z)$, and acts on $K_*(\pi_{ku})
\otimes \Q \cong \Q[u] \{ du \}$ by sign on $u$ and $du$, hence fixes
$\Q[u^2] \{u du\} \cong \Q[\alpha] \{d\alpha\} \cong K_*(\pi_{ko})
\otimes \Q$.  So $K(ko) \to K(ku)^{h\Z/2}$ is a rational equivalence.
The conjugation action also fixes $BK(\Z)$ after rationalization, so the
Poincar{\'e} series of $K(KU)^{h\Z/2}$ is $(1+t) + (2t^5 + t^6)/(1-t^4)$.
\endremark

\remark{Remark 2.14}
We expect that $c_*$ and $\pi_*$ are essential (not null-homotopic) as
maps of $A(S^3)$-module spectra and $K(ku)$-module spectra, respectively.
In other words, we expect that $K(ku) \to K(KU) \to \Sigma K(\Z)$ is
a non-split extension of $K(ku)$-module spectra.  This expectation is
to some extent justified by the fact that the cofiber $THH(ku|KU)$
of the $THH$-transfer map $\pi_* \: THH(\Z) \to THH(ku)$ sits in
a non-split extension $THH(ku) \to THH(ku|KU) \to \Sigma THH(\Z)$
of $THH(ku)$-module spectra.  See \cite{Au05, 10.4}, or \cite{HM03,
Lemma~2.3.3} for a similar result in an algebraic case.
\endremark

\head \S3. Examples from smooth bordism \endhead

\subhead 3.1. Oriented bordism spectra \endsubhead
Let $MU$ be the complex bordism spectrum, $MSO$ the real oriented
bordism spectrum, and $MSp$ the symplectic bordism spectrum.  These are
all connective commutative $S$-algebras, given by the Thom
spectra associated to infinite loop maps from $BU$, $BSO$ and $BSp$
to $BSF = BSL_1(S)$, respectively.  We recall that
$$
H_*(BU) \cong \Z[b_k \mid k\ge1]
$$
with $|b_k|=2k$, while $H_*(BSO; \Z[1/2]) \cong H_*(BSp; \Z[1/2]) \cong
\Z[1/2][q_k \mid k\ge1]$ with $|q_k| = 4k$.

The Thom equivalence $\theta \: MU \wedge MU \to MU \wedge S[BU]$
induces an equivalence $H\Z \wedge MU \simeq H\Z \wedge S[BU] = H\Z[BU]$.
Combined with the Hurewicz map $\pi \: S \to H\Z$ we obtain a chain of
maps of commutative $S$-algebras
$$
MU @>>> H\Z \wedge MU \simeq H\Z [BU] @<<< S[BU] \,,
$$
that are $\pi_0$-isomorphisms and rational equivalences.  There are
similar chains $MSO @>>> H\Z \wedge MSO \simeq H\Z [BSO] @<<< S[BSO]$
and $MSp @>>> H\Z \wedge MSp \simeq H\Z [BSp] @<<< S[BSp]$, and all
induce rational equivalences
$$
\aligned
K(MU) &@>>> K(\Z[BU]) @<<< A(BBU) \\
K(MSO) &@>>> K(\Z[BSO]) @<<< A(BBSO) \\
K(MSp) &@>>> K(\Z[BSp]) @<<< A(BBSp)
\endaligned
\tag 3.2
$$
of commutative $S$-algebras.  Here we view $BU \simeq \Omega BBU$ as
the Kan loop group of $BBU$, $\Z[BU]$ is the associated simplicial ring,
and similarly for $BSO$ and $BSp$.

\proclaim{Lemma 3.3}
(a)
The simplicial $\Q$-algebra $R = \Q[BU]$ with $\pi_* R = H_*(BU; \Q) =
\Q[b_k \mid k\ge1]$ has Hochschild homology
$$
HH_*(R) \cong \Q[b_k \mid k\ge1] \otimes E(db_k \mid k\ge1) \,,
$$
with Poincar{\'e} series
$$
h(t) = \prod_{k\ge1} {1 + t^{2k+1} \over 1 - t^{2k}} \,,
$$
and Connes' operator acts by $B(b_k) = db_k$.

(b) The de\,Rham homology $H^{dR}_*(R) \cong \Q$ is concentrated in
degree zero, so $\pi_R \: R \to \Q$ is a de\,Rham equivalence.

(c) The Poincar{\'e} series of
$\ker(B) \subset HH_*(R)$ is
$$
k(t) = {1 + t h(t) \over 1 + t} \,.
$$

(d)
The simplicial $\Q$-algebra $R_{so} = \Q[BSO] \simeq \Q[BSp]$,
with $\pi_* R_{so} = \Q[q_k \mid k\ge1]$, has Hochschild homology
$$
HH_*(R_{so}) \cong \Q[q_k \mid k\ge1] \otimes E(dq_k \mid k\ge1)
\,.
$$
Its Poincar{\'e} series is $h_{so}(t) = \prod_{k\ge1} (1 + t^{4k+1}) /
(1 - t^{4k})$.  The map $R_{so} \to \Q$ is a de\,Rham equivalence, and
$\ker(B) \subset HH_*(R_{so})$ has Poincar{\'e} series $k_{so}(t) =
(1 + t h_{so}(t)) / ( 1+t )$.
\endproclaim

\demo{Proof}
(a) In this case the spectral sequence~(2.4) has $E^2_{**} =
HH_*(\Q[b_k \mid k\ge1]) \cong \Q[b_k \mid k\ge1] \otimes E(db_k \mid
k\ge1)$.  The algebra generators are in filtrations $0$ and $1$, so $E^2
= E^\infty$.  This term is free as a graded commutative $\Q$-algebra,
so $HH_*(R)$ is isomorphic to the $E^\infty$-term.

(b) The homology of $\Q[b_k] \otimes E(db_k)$ with respect to $B$ is just
$\Q$, for each $k\ge1$, so by the K{\"u}nneth theorem the de\,Rham homology
of $HH_*(R)$ is also just $\Q$.

(c) Write $H_n$ for $HH_n(R)$ and $K_n$ for $\ker(B \: H_n \to H_{n+1})$.
Let $h_n = \dim_\Q H_n$, so $h(t) = \sum_{n\ge0} h_n t^n$, and $k_n =
\dim_\Q K_n$.  In view of the exact sequence
$$
0 \to \Q \to H_0 @>d>> H_1 @>d>> \dots @>d>> H_{n-1} \to K_n \to 0
$$
we find that $1 - (-1)^n k_n = h_0 - h_1 + \dots + (-1)^{n-1} h_{n-1}$, so
$$
\sum_{n\ge0} t^n k_n - \sum_{n\ge0} (-t)^n
= t h(t) - t^2 h(t) + \dots + (-t)^{m+1} h(t) + \dots \,.
$$
It follows that the Poincar{\'e} series $k(t) = \sum_{n\ge0} k_n t^n$
for $\ker(B)$ satisfies $k(t) - 1/(1+t) = t h(t)/(1+t)$.

(d) The only change from the complex to the oriented real and symplectic
cases is in the grading of the algebra generators, which (as long as
they remain in even degrees) plays no role for the proofs.  \qed
\enddemo

\proclaim{Theorem 3.4}
(a) There is a rationally split homotopy fiber sequence
$$
K(\pi_{MU}) \to K(MU) @>\pi>> K(\Z)
$$
and the trace map $tr \: K(MU) \to THH(MU)$ identifies
$K_*(\pi_{MU}) \otimes \Q$ with the positive-degree part of $\ker(B)$
in
$$
THH_*(MU) \otimes \Q \cong \Q[b_k \mid k\ge1] \otimes E(db_k \mid k\ge1)
\,,
$$
where $|b_k| = 2k$ and $B(b_k) = db_k$.  Hence $K(\pi_{MU})$
has Poincar{\'e} series
$$
k(t)-1 = {t h(t)-t \over 1+t} \,.
$$

(b) There are rationally split homotopy fiber sequences
$$
\align
K(\pi_{MSO}) &\to K(MSO) @>\pi>> K(\Z) \\
K(\pi_{MSp}) &\to K(MSp) @>\pi>> K(\Z)
\endalign
$$
and the trace maps identify both $K_*(\pi_{MSO}) \otimes \Q$ and
$K_*(\pi_{MSp}) \otimes \Q$ with the positive-degree part of $\ker(B)$ in
$$
THH_*(MSO) \otimes \Q \cong THH_*(MSp) \otimes \Q
\cong \Q[q_k \mid k\ge1] \otimes E(dq_k \mid k\ge1) \,,
$$
where $|q_k| = 4k$ and $B(q_k) = dq_k$.  Hence $K(\pi_{MSO})$ and
$K(\pi_{MSp})$ both have Poincar{\'e} series
$k_{so}(t)-1 = (t h_{so}(t)-t )/( 1+t )$.
\endproclaim

\remark{Remark 3.5}
Adding the Poincar{\'e} series of $K(\Z)$, as in Remark~2.6, we find
that the Poincar{\'e} series of $K(MU)$ and $A(BBU)$ is
$$
{t^5 \over 1-t^4} + {1 + t h(t) \over 1+t} \,,
$$
whereas the Poincar{\'e} series of $K(MSO)$, $A(BBSO)$,
$K(MSp)$ and $A(BBSp)$ is $t^5/(1-t^4) + (1 + t h_{so}(t))/(1+t)$.
\endremark

\head \S4. Units, determinants and traces \endhead

\subhead 4.1. Units \endsubhead
For each connective $S$-algebra $A$ there is a natural map of spaces
$$
w \: BGL_1(A) \to K(A)
$$
that factors as the infinite stabilization map $BGL_1(A) \to
BGL_\infty(A)$, composed with the inclusion $BGL_\infty(A) \to
BGL_\infty(A)^+$ into Quillen's plus construction, and followed by the
inclusion of $BGL_\infty(A)^+ \cong \{1\} \times BGL_\infty(A)^+$ into
$K_0(\pi_0 A) \times BGL_\infty(A)^+ = K(A)$.

\remark{Remark 4.2}
This $w$ is an $E_\infty$ map with respect to the multiplicative
$E_\infty$ structure on $K(A)$ that is induced by the smash product over
$A$.  However, we shall only work with the additive grouplike $E_\infty$
structure on $K(A)$, which comes from viewing $K(A)$ as the underlying
infinite loop space of the $K$-theory spectrum.  So when we refer to
infinite loop structures below, we are thinking of the additive ones.
\endremark

\medskip

We write $BSL_1(A) = BGL_1(\pi_A)$ for the homotopy fiber of the map
$BGL_1(A) \to BGL_1(\pi_0 A)$ induced by $\pi_A \: A \to H\pi_0 A$.
In the resulting diagram
$$
BSL_1(A) @>w>> K(A) @>\pi>> K(\pi_0 A)
\tag 4.3
$$
the composite map has a preferred null-homotopy (to the base point
of the $1$-component of $K(\pi_0 A)$).  The diagram is a rational
homotopy fiber sequence if and only if $w \: BSL_1(A) \to K(\pi_A)$
is a rational equivalence.  Note that the natural inclusion $\{1\}
\times BGL_\infty(A)^+ \to K(A)$ induces a homotopy equivalence
$$
BGL_\infty(\pi_A)^+ \simeq K(\pi_A) \,,
$$
since $K_0(A) \cong K_0(\pi_0 A)$.

\subhead 4.4. Determinants \endsubhead
Suppose furthermore that $A$ is commutative as an $S$-algebra.  One
attempt at proving that $w$ is injective could be to construct a map $\det
\: K(A) \to BGL_1(A)$ with the property that $\det \circ w \simeq \id$.
However, no such determinant map exists in general, as the following
adaption of \cite{Wa82, 3.7} shows.

\example{Example 4.5}
When $A = S$, the map $\lambda \circ w \: BF = BGL_1(S) \to A(*) \to
K(\Z)$ factors through $BGL_1(\Z) \simeq K(\Z/2, 1)$, and $\pi_2(\lambda)
\: \pi_2 A(*) \to K_2(\Z) \cong \Z/2$ is an isomorphism, so $\pi_2(w) \:
\pi_2 BF \to \pi_2 A(*)$ is the zero map.  But $\pi_2 BF \cong \pi_1(S)
\cong \Z/2$ is not zero, so $\pi_2(w)$ is not injective.  In particular,
$w$ is not split injective up to homotopy.
\endexample

However, it is possible to construct a rationalized determinant map.
Recall from Subsection~1.1 that $A_\Q$ is equivalent to $HR$ for
some naturally determined commutative simplicial $\Q$-algebra $R$.

\proclaim{Lemma 4.6}
Let $R$ be a commutative simplicial ring.  There is a natural infinite
loop map
$$
\det \: BGL_\infty(R)^+ \to BGL_1(R)
$$
that agrees with the usual determinant map for discrete commutative
rings, such that the composite with $w \: BGL_1(R) \to BGL_\infty(R)^+$
equals the identity.
\endproclaim

\demo{Proof}
The usual
matrix determinant $\det \: M_n(R) \to R$ induces a simplicial
group homomorphism $GL_n(R) \to GL_1(R)$ and a pointed map
$BGL_n(R) \to BGL_1(R)$ for each $n\ge0$.  These stabilize to a
map $BGL_\infty(R) \to BGL_1(R)$, which extends to an
infinite loop map
$$
\det \: BGL_\infty(R)^+ \to BGL_1(R) \,,
$$
unique up to homotopy, by the multiplicative infinite loop structure
on the target and the universal property of Quillen's plus construction.
To make the construction natural, we fix a choice of extension
in the initial case $R = \Z$, and define $\det_R$ for general $R$
as the dashed pushout map in the following diagram:
$$
\xymatrix{
BGL_\infty(\Z) \ar[r] \ar[d] & BGL_\infty(\Z)^+ \ar[d] \ar[dr]^{\det_\Z} \\
BGL_\infty(R) \ar[r] \ar[drr] & BGL_\infty(R)^+ \ar@{-->}[dr]^{\det_R} 
	& BGL_1(\Z) \ar[d] \\
& & BGL_1(R)
}
$$
(Recall that Quillen's plus construction is made functorial by demanding
that the left hand square is a pushout.)
\qed
\enddemo

\proclaim{Proposition 4.7}
Let $A$ be a connective commutative $S$-algebra.  There is a natural
infinite loop map
$$
{\det}_\Q \: BGL_\infty(A)^+ \to BGL_1(A)_\Q
$$
that agrees with the rationalized determinant map for a commutative ring
$R$ when $A = HR$, such that the composite
$$
BGL_1(A) @>w>> BGL_\infty(A)^+ @>\det_\Q>> BGL_1(A)_\Q
$$
is homotopic to the rationalization map.
\endproclaim

\demo{Proof}
We define $\det_\Q$ as the dashed pullback map in the following diagram
$$
\xymatrix{
BGL_\infty(A)^+ \ar[r] \ar[d] \ar@{-->}[dr]^{\det_\Q} & BGL_\infty(A_\Q)^+
	\ar[dr]^{\det'_R} \\
BGL_\infty(\pi_0 A)^+ \ar[dr]_{(\det_{\pi_0 A})_\Q}
	& BGL_1(A)_\Q \ar[r] \ar[d] & BGL_1(A_\Q) \ar[d] \\
& BGL_1(\pi_0 A)_\Q \ar[r] & BGL_1(\pi_0 A_\Q)
}
$$
where the vertical maps are induced by the Postnikov section $\pi \: A
\to H\pi_0 A$, and the horizontal maps are induced by the rationalization
$q \: A \to A_\Q$.  The right hand square is a homotopy pullback, since
$SL_1(A)_\Q \simeq SL_1(A_\Q)$.

To define the map $\det'_R$, we take $R$ to be a commutative simplicial
$\Q$-algebra such that $A_\Q \simeq HR$ as commutative $H\Q$-algebras.
A natural choice can be made for $R$, as discussed in Subsection~1.1,
such that the identification $\pi_0 A_\Q \cong \pi_0 R$ is the identity.
Then $\det'_R$ is the composite map
$$
BGL_\infty(A_\Q)^+ \simeq BGL_\infty(R)^+
@>\det_R>> BGL_1(R) \simeq BGL_1(A_\Q) \,,
$$
with $\det_R$ from Lemma~4.6.  It strictly covers the map
$\det_{\pi_0 A_\Q}$, so the outer hexagon commutes strictly.
This defines the desired map $\det_\Q$.

To compare $\det_\Q \circ w$ and $q \: BGL_1(A) \to BGL_1(A)_\Q$,
note that both maps have the same composite to $BGL_1(\pi_0 A)_\Q$,
they have homotopic composites to $BGL_1(A_\Q)$, and all composites
(and homotopies) to $BGL_1(\pi_0 A_\Q)$ are equal.  Hence the maps to
the homotopy pullback are homotopic, too.
\qed
\enddemo

\proclaim{Theorem 4.8}
(a) The relative unit map
$$
BBU_\otimes = BGL_1(\pi_{ku})
	@>w>> BGL_\infty(\pi_{ku})^+ \simeq K(\pi_{ku})
$$
is a rational equivalence, with rational homotopy inverse given by the
relative rational determinant map
$$
{\det}_\Q \: BGL_\infty(\pi_{ku})^+
	\to BGL_1(\pi_{ku})_\Q = (BBU_\otimes)_\Q \,.
$$

(b) The relative unit maps
$$
\align
BBO_\otimes = BGL_1(\pi_{ko}) &\to K(\pi_{ko}) \\
BW_\otimes = BGL_1(\pi_\ell) &\to K(\pi_\ell)
\endalign
$$
are rational equivalences (with rational homotopy inverse $\det_\Q$
in each case).
\endproclaim

\demo{Proof}
(a)
By Proposition~4.7, the composite
$$
BBU_\otimes @>w>> K(\pi_{ku}) @>\det_\Q>> (BBU_\otimes)_\Q
$$
is a rational equivalence, so $w$ is rationally injective.  Here $\pi_*
BBU_\otimes \cong \pi_{*-1} BU_\otimes$ has Poincar{\'e} series
$t^3/(1-t^2)$, just like $K(\pi_{ku})$ by Theorem~2.5(a).  Thus $w$
is a rational equivalence.

(b)
The same proof works for $ko$ and $\ell$, using that $BBO_\otimes$ and
$BW_\otimes$ have Poincar{\'e} series $t^5/(1-t^4)$ and
$t^{2p-1}/(1-t^{2p-2})$, respectively.
\qed
\enddemo

\remark{Remark 4.9}
The analogous map $w \: BSL_1(MU) \to K(\pi_{MU})$ is rationally
injective, but not a rational equivalence.  For the Poincar{\'e} series
of the source is
$$
t(p(t)-1) = t^3 + 2t^5 + 3t^7 + 5t^9 + \dots \,,
$$
where $p(t) = \prod_{k\ge1} 1/(1-t^{2k})$,
and the Poincar{\'e} series of the target is
$$
(th(t)-t)/(1+t) = t^3 + 2t^5 + 3t^7 + t^8 + 5t^9 + \dots \,,
$$
by Theorem~3.4(a).  These first differ in degree~8, since $\pi_8
BSL_1(MU) \cong \pi_7 MU$ is trivial, but $K_8(MU)$ and $K_8(\pi_{MU})$ 
have rank one.  A generator of the latter group maps to $db_1 \cdot db_2$
in $\ker(B) \subset THH_*(MU) \otimes \Q$.

In the same way, $w \: BSL_1(MSO) \to K(\pi_{MSO})$ and its symplectic
variant are rationally injective, but not rational equivalences.
\endremark

\subhead 4.10. Traces \endsubhead
Our original strategy for proving that $w \: BGL_1(A) \to K(A)$ is
rationally injective for $A = ku$ was to use the trace map $tr \: K(A)
\to THH(A)$, in place of the rational determinant map.
By \cite{Schl04, \S4}, there is a natural commutative diagram
$$
\xymatrix{
BGL_1(A) \ar[r]^w \ar[d] & K(A) \ar[d] \ar[dr]^{tr} \\
B^{cy}GL_1(A) \ar[r] & K^{cy}(A) \ar[r] & THH(A) \\
GL_1(A) \ar[rr] \ar[u] & & \Omega^\infty A \ar[u]
}
$$
where $B^{cy}$ and $K^{cy}$ denote the cyclic bar construction and cyclic
$K$-theory, respectively.  The middle row is the geometric realization
of two cyclic maps, hence consists of circle equivariant spaces and maps.

When $A = ku$, the resulting $B$-operator on $H_*(B^{cy} BGL_1(A); \Q)$
takes primitive classes in the image from $H_*(GL_1(A); \Q) \cong
H_*(BU_\otimes; \Q)$ to primitive classes generating the image from
$H_*(BGL_1(A); \Q) \cong H_*(BBU_\otimes; \Q)$, so by a diagram chase we
can determine the images of the latter primitive classes in $H_*(THH(A);
\Q)$.  By an appeal to the Milnor--Moore theorem \cite{MM65, App.}, this
suffices to prove that $tr \circ w$ is rationally injective in this case.

In comparison with the rational determinant approach taken above,
this trace method involves more complicated calculations.  For
commutative $S$-algebras, it is therefore less attractive.  However,
for non-commutative $S$-algebras, the trace method may still be useful,
since no (rational) determinant map is likely to exist.  We have therefore
sketched the idea here, with a view to future applications.

\head \S5. Two-vector bundles and elliptic objects \endhead

The following discussion elaborates on the second author's work with
Baas and Dundas in \cite{BDR04}.  It is intended to explain some
of our interest in Theorem~0.1.

\subhead 5.1. Two-vector bundles \endsubhead
A $2$-vector bundle $\Cal E$ of rank~$n$ over a base space $X$ is
represented by a map $X \to |BGL_n(\Cal V)|$, where $\Cal V$ is the
symmetric bimonoidal category of finite dimensional complex vector
spaces.  A virtual $2$-vector bundle $\Cal E$ over $X$ is represented
by a map $X \to K(\Cal V)$, where $K(\Cal V)$ the algebraic $K$-theory
of the $2$-category of finitely generated free $\Cal V$-modules; see
\cite{BDR04, Thm.~4.10}.  By \cite{BDRR}, spectrification induces a
weak equivalence $\Spt \: K(\Cal V) \to K(ku)$, so the $2$-vector bundles
over $X$ are geometric $0$-cycles for the cohomology theory $K(ku)^*(X)$.

\subhead 5.2. Anomaly bundles \endsubhead
The preferred rational splitting of $\pi \: K(ku) \to K(\Z)$ defines an
infinite loop map
$$
{\det}_\Q \: K(ku) \to (BBU_\otimes)_\Q \,,
$$
which extends the rationalization map over $w \: BBU_\otimes \to K(ku)$
and agrees with the relative rational determinant on $K(\pi_{ku})$.  (We
do not know if there exists an integral determinant map $BGL_\infty(ku)^+
\to BBU_\otimes$ in this case.)  We define the rational determinant bundle
$|\Cal E| = \det(\Cal E)$ of a virtual $2$-vector bundle represented by
a map $\Cal E \: X \to K(\Cal V) \simeq K(ku)$, as the composite map
$$
|\Cal E| \: X @>\Cal E>> K(ku) @>\det_\Q>> (BBU_\otimes)_\Q \,.
$$
We define the rational anomaly bundle $\Cal H \to \Cal L X$ of $\Cal E$
as the composite map
$$
\Cal H \: \Cal L X @>\Cal L |\Cal E|>> \Cal L (BBU_\otimes)_\Q
	@>r_\Q>> (BU_\otimes)_\Q \,,
$$
where $r \: \Cal L BBU_\otimes \to BU_\otimes$ is the retraction defined
as the infinite loop cofiber of the constant loops map $BBU_\otimes
\to \Cal L BBU_\otimes$.  Up to rationalization, $\Cal H$ is a virtual
vector bundle of virtual dimension $+1$, i.e., a virtual line bundle.
Furthermore, the anomaly bundle relates the composition $\star$ of free
loops, when defined, to the tensor product of virtual vector spaces:
the square
$$
\xymatrix{
\Cal L X \times_X \Cal L X \ar[r]^-{(\Cal H, \Cal H)} \ar[d]_\star &
(BU_\otimes)_\Q \times (BU_\otimes)_\Q \ar[d]^\otimes \\
\Cal L X \ar[r]^{\Cal H} & (BU_\otimes)_\Q
}
\tag 5.3
$$
commutes up to coherent isomorphism.

\subhead 5.4. Gerbes \endsubhead
A $2$-vector bundle of rank~$1$ over $X$ is the same as a $\C^*$-gerbe
$\Cal G$, which is represented by a map $\Cal G \: X \to BBU(1)$.  When viewed
as a virtual $2$-vector bundle, via $BBU(1) \to BBU_\otimes
\to K(ku)$, the associated anomaly bundle is the complex
line bundle over $\Cal L X$ that is represented by the composite
$$
\Cal L X @>\Cal L \Cal G>> \Cal L BBU(1) @>r>> BU(1) \,.
$$
This is precisely the anomaly line bundle for $\Cal G$, as described in
\cite{Br93, \S6.2}.  Note that the rational anomaly bundles
of virtual 2-vector bundles represent
general elements in
$$
1 + \widetilde K^0(\Cal L X) \otimes \Q
	\subset K^0(\Cal L X) \otimes \Q \,,
$$
whereas the anomaly line bundles of gerbes only represent 
elements in $H^2(\Cal L X)$.

\subhead 5.5. State spaces and action functionals \endsubhead
In physical language, we think of a free loop $\gamma \: S^1 \to X$
as a closed string in a space-time $X$.  For a $2$-vector bundle $\Cal
E \to X$, we think of the fiber $H_\gamma$ (a virtual vector space) at
$\gamma$ of the anomaly bundle $\Cal H \to \Cal L X$ as the state space
of that string.  Then the state space of a composite of two strings
(or a disjoint union of two strings) is the tensor product of the
individual state spaces, as is usual in quantum mechanics.  Similarly,
the state space of an empty set of strings is $\C$.  In the special case
of anomaly line bundles for gerbes, the resulting state spaces are only
complex lines, but in our generality they are virtual vector spaces.
These are much closer to the Hilbert spaces usually considered in
more analytical approaches to this subject.

There is evidence that a two-part differential-geometric structure
$(\nabla_1, \nabla_2)$ on $\Cal E$ over $X$ (somewhat like a
connection for a vector bundle, but providing parallel transport both
for objects and for morphisms in the $2$-vector bundle) provides $\Cal
H \to \Cal L X$ with a connection, and more generally an action functional
$$
S(\Sigma) \: H_{\bar\gamma_1} \otimes \dots \otimes H_{\bar\gamma_p}
\to H_{\gamma_1} \otimes \dots \otimes H_{\gamma_q} \,,
$$
where $\Sigma \: F \to X$ is a
compact Riemann surface over $X$, with $p$ incoming and $q$ outgoing
boundary circles.  The time development of the physical system is then
given by the Euler--Lagrange equations of the action functional.

In a little more detail, the idea is that the primary form of parallel
transport in $(\Cal E, \nabla_1)$ around $\gamma$ provides an endo-functor
$\tilde\gamma$ of the fiber category $\Cal E_x \cong \Cal V^n$ over a
chosen point $x$ of $\gamma$.  More precisely, parallel transport only
provides a zig-zag of functors connecting $\Cal E_x$ to itself, but the
determinant in $(BU_\otimes)_\Q$ is still well-defined.  This ``holonomy''
is then the fiber $H_\gamma = \det(\tilde\gamma)$ at $\gamma$ of the
anomaly bundle $\Cal H$.  For a moving string, say on the Riemann surface
$F$, the secondary form of parallel transport $\nabla_2$ specifies how
the holonomy changes with the string, and this defines the connection
$\nabla$ on $\Cal H \to \Cal L X$.  In the gerbe case, this theory has
been worked out in \cite{Br93, \S5.3}, where $\nabla_1$ is called
``connective structure'' and $\nabla_2$ is called ``curving''.

For a closed surface $F$, $S(\Sigma) \: \C \to \C$ is multiplication by
a complex number, which would only depend on the rational type
of~$\Cal E$.  Optimistically, this association can produce a conformal
invariant of $F$ over $X$, which in the case of genus~$1$ surfaces would
lead to an elliptic modular form.  Less naively, additional structure
derived from a string structure on $X$ should account for the weight of
the modular form.  With such structure, a $2$-vector bundle $\Cal E$
with connective structure $\nabla$ would qualify as a Segal elliptic
object over $X$.

\subhead 5.6. Open strings \endsubhead
In the presence of $D$-branes in the space-time $X$, we can extend the
anomaly bundle to also cover open strings with end points restricted
to lie on these $D$-branes; see \cite{Mo04, \S3.4}.  In this terminology,
the (rational) determinant bundle $|\Cal E| \to X$ plays the
role of the $B$-field.

By a (rational) {\bf $D$-brane} $(\Cal W, E)$ in $X$ we will mean a
subspace $\Cal W \subset X$ together with a trivialization $E$ of the
restriction of the (rational) determinant bundle $|\Cal E|$ to $\Cal W$.
In terms of representing maps, $E$ is a null-homotopy of the composite map
$$
\Cal W \subset X @>\Cal E>> K(ku) @>\det_\Q>> (BBU_\otimes)_\Q \,.
$$
In similar terminology, we may refer to the determinant bundle $|\Cal E|
\to X$ as the (rational) {\bf $B$-field}.

When the $B$-field $|\Cal E|$ is rationally trivial, then a second
choice of trivialization $E$ amounts to a choice of null-homotopy of the
trivial map $\Cal W \to (BBU_\otimes)_\Q$, or equivalently to a map $E \:
\Cal W \to (BU_\otimes)_\Q$.  In other words, $E$ is a virtual vector
bundle over $\Cal W$ of virtual dimension $+1$, up to rationalization.
In this case, the $K$-theory class of $E \to \Cal W$ in $1 + \widetilde
K^0(\Cal W) \otimes \Q$ is the ``charge'' of the $D$-brane $(\Cal W, E)$.
This conforms with the (early) view on $D$-branes as coming equipped
with a charge $[E]$ in topological $K$-theory.

For a general $B$-field $|\Cal E|$, the possible trivializations $E$ of
its restriction to $\Cal W$ instead form a torsor under the group $1 +
\widetilde K^0(\Cal W) \otimes \Q$.  For two such trivializations $E$
and $E'$ differ by a loop of maps $\Cal W \to (BBU_\otimes)_\Q$, or
equivalently a map $E' - E \: \Cal W \to (BU_\otimes)_\Q$.  So $[E' -
E]$ is a topological $K$-theory class measuring the charge difference
between the two $D$-branes $(\Cal W, E)$ and $(\Cal W, E')$.

Given two $D$-branes $(\Cal W_0, E_0)$ and $(\Cal W_1, E_1)$ in
$(X, \Cal E)$, we have a commutative diagram
$$
\xymatrix{
\Cal W_0 \ar[r] \ar[d]_{E_0} & X \ar[d]^{|\Cal E|}
  & \Cal W_1 \ar[l] \ar[d]^{E_1} \\
P (BBU_\otimes)_\Q \ar[r]^{\pi} & (BBU_\otimes)_\Q
  & P (BBU_\otimes)_\Q \ar[l]_{\pi}
}
$$
where $\pi \: PY \to Y$ denotes the path space fibration covering a
based space $Y$.  An open string in $X$, constrained to $\Cal W_0$ and
$\Cal W_1$ at its ends, is a map $\gamma \: I \to X$ with $\gamma(0)
\in \Cal W_0$ and $\gamma(1) \in \Cal W_1$.  In other words, it is an
element in the homotopy pullback of the top row in the diagram above.
Let $\Omega(X, \Cal W_0, \Cal W_1)$ denote the space of such open strings.
The homotopy pullback of the lower row is $\Omega (BBU_\otimes)_\Q \simeq
(BU_\otimes)_\Q$.  Hence the $2$-vector bundle $\Cal E$ and the two $D$-branes
specify a map of homotopy pullbacks
$$
\Cal H \: \Omega(X, \Cal W_0, \Cal W_1) \to (BU_\otimes)_\Q
$$
that we call the (rational, virtual) anomaly bundle of this space of
open strings.  Again, we think of each fiber $H_\gamma$ at $\gamma \: (I,
0, 1) \to (X, \Cal W_0, \Cal W_1)$ as the state space of that open string.

In the presence of a suitable connection $(\nabla_1, \nabla_2)$ on $\Cal
E \to X$, parallel transport in $(\Cal E, \nabla_1)$ along $\gamma$
induces a (zig-zag) functor $\tilde\gamma$ from $\Cal E_x$
to $\Cal E_y$, with determinant $\det(\tilde\gamma)$
from the fiber of $|\Cal E|$ at $x=\gamma(0)$ to the fiber at
$y=\gamma(1)$.  The trivializations of these two
fibers provided by the $D$-brane data $E_0$ and $E_1$, respectively,
then agree up to a correction term, which is the fiber $H_\gamma$
in the anomaly bundle:
$$
\det(\tilde\gamma)(E_{0,x}) \cong H_\gamma \otimes E_{1,y}
$$
Again, the secondary part of the connection may induce a connection on
$\Cal H$ over $\Omega(X, \Cal W_0, \Cal W_1)$, and more generally an
action functional $S(\Sigma)$, where now $\Sigma \: F \to X$ and the
incoming and the outgoing parts of $F$ are unions of circles and closed
intervals.  For example, an open string might split off a closed string.
One advantage of the above perspective is that the state spaces of open
and closed strings arise in a compatible fashion, as the holonomy of
parallel transport in the $2$-vector bundle $\Cal E$, and this makes
the construction of $S(\Sigma)$ feasible.  The gerbe case is discussed
in \cite{Br93, \S6.6}.

\enddocument